\documentclass[12pt]{amsart}
\usepackage[utf8]{inputenc}
\usepackage[english]{babel}
\usepackage{color}
\usepackage{amssymb}
\usepackage{amsmath}
\usepackage{amsthm}
\usepackage{leftidx}
\usepackage{mathdots}
\usepackage{amssymb}
\usepackage{mathrsfs}

\usepackage{epsfig}
\usepackage{verbatim}
\newcommand{\lra}{{\longrightarrow}}
\newcommand{\eproof}{\hfill\rule{2.2mm}{3.0mm}}

\newcommand{\Proof}{\noindent {\bf Proof.~~}}

\newcommand{\R}{{\mathbb R}}
\newcommand{\Z}{{\mathbb Z}}

\newcommand{\C}{{\mathbb C}}

\newcommand{\PP}{{\mathbb P}}

\renewcommand{\eqref}[1]{(\ref{#1})}
\newcommand{\inner}[1]{\langle #1 \rangle}
\newcommand{\innerp}[1]{\langle #1 \rangle}
\newcommand{\shsp}{\hspace{1em}}
\newcommand{\mhsp}{\hspace{2em}}
\newcommand{\abs}[1]{\lvert#1\rvert}
\newcommand{\AT}[1]{\widehat{#1}}

\newcommand{\A}{{\mathcal A}}

\newcommand{\rank}{{\rm rank}}
\newcommand{\Span}{{\rm span}}

\newcommand{\tr}{{\rm tr}}

\newcommand{\vx}{{\mathbf x}}
\newcommand{\vy}{{\mathbf y}}
\newcommand{\vz}{{\mathbf z}}
\newcommand{\vw}{{\mathbf w}}

\newcommand{\vv}{{\mathbf v}}
\newcommand{\vu}{{\mathbf u}}
\newcommand{\vf}{{\mathbf f}}

\newcommand{\ve}{{\mathbf e}}
\renewcommand{\H}{{\mathbb F}}
\newcommand{\HH}{{\H}^d}
\newcommand{\qHH}{\underline{\HH}}
\newcommand{\ul}[1]{\underline{#1}}
\newcommand{\MM}{\mathbf M}

\newcommand{\Herm}{{\mathbf H}_d}

\newcommand{\uD}{\underline{D}}
\newtheorem{prop}{Proposition}[section]
\newtheorem{lem}[prop]{Lemma}
\newtheorem{defi}{Definition}[section]
\newtheorem{coro}[prop]{Corollary}
\newtheorem{theo}[prop]{Theorem}
\newtheorem{remark}[prop]{Remark}

\begin{document}
\baselineskip 18pt
\title[Almost Everywhere Generalized Phase Retrieval]{Almost Everywhere Generalized Phase Retrieval}

\author{Meng Huang}
\address{LSEC, Inst.~Comp.~Math., Academy of
Mathematics and System Science,  Chinese Academy of Sciences, Beijing, 100190, China}
\email{hm@lsec.cc.ac.cn}

\author{Yi Rong}
\address{Department of Mathematics  \\ Hong Kong University of Science and Technology\\
Clear Water Bay, Kowloon, Hong Kong}
\email{yrong@ust.hk}

\author{Yang Wang}
\thanks{Yang Wang was supported in part by the Hong Kong Research Grant Council grant 16306415.
      Zhiqiang Xu was supported  by NSFC grant (91630203, 11688101),
Beijing Natural Science Foundation (Z180002).}
\address{Department of Mathematics  \\ Hong Kong University of Science and Technology\\
Clear Water Bay, Kowloon, Hong Kong}
\email{yangwang@ust.hk}
\author{Zhiqiang Xu}
\address{LSEC, Inst.~Comp.~Math., Academy of
Mathematics and System Science,  Chinese Academy of Sciences, Beijing, 100190, China\newline
School of Mathematical Sciences, University of Chinese Academy of Sciences, Beijing 100049, China}
\email{xuzq@lsec.cc.ac.cn}

\begin{abstract}
    The aim of generalized phase retrieval is to recover $\vx\in \H^d$ from the quadratic measurements $\vx^*A_1\vx,\ldots,\vx^*A_N\vx$,  where $A_j\in \Herm(\H)$ and $\H=\R$ or $\C$. In this paper, we study the matrix set $\A=(A_j)_{j=1}^N$ which has
    the  almost everywhere phase retrieval property. For the case $\H=\R$, we show that $N\geq d+1$ generic matrices with   prescribed ranks have almost everywhere phase retrieval property. We also extend this result to the case where $A_1,\ldots,A_N$ are orthogonal matrices and hence establish the almost everywhere phase retrieval property for the fusion frame phase retrieval.   For the case where $\H=\C$, we obtain similar results under  the assumption  of  $N\geq 2d$.
 We lower the measurement number $d+1$ (resp. $2d$) with showing that there exist $N=d$ (resp. $2d-1$) matrices $A_1,\ldots, A_N\in \Herm(\R)$ (resp. $\Herm(\C)$) which have the almost everywhere phase retrieval property.
 Our results are an extension of almost everywhere phase retrieval  from the standard phase retrieval to the general setting and the proofs are often based on some new ideas about determinant variety.
\end{abstract}
\maketitle

\section{Introduction}
\setcounter{equation}{0}

\subsection{Problem Setup}
The classic phase retrieval problem, which concerns the reconstruction of a function from the magnitude of its Fourier transform, has many applications in various areas such as X-ray crystallography, super-resolution cryo-EM imaging, optics, signal processing and many more.  It is well known that the map $f \mapsto \AT f$, where $\AT f$ denotes the Fourier transform of $f$, is an isometry in $L^2(\R^d)$ and hence $f$ can be uniquely reconstructed from $\AT f$. However, when only the magnitude $|\AT f(\xi)|$ is known, the reconstruction becomes rather nontrivial. In many cases, the reconstruction is impossible because the solution is not unique, even after taking some obvious factors into consideration, such as translation and modulation.

More recently the phase retrieval problem has been naturally extended to finite dimensional Hilbert spaces, and research in this area have accounted for the bulk of the advances lately (see e.g. \cite{BCE06,bodmann,CSV12,CEHV15, FMW14, HMW13} and the references therein).
We mainly focus on the finite dimensional Hilbert space $\H^d$ where $\H=\R$ or $\H=\C$ and consider a set of Hermitian matrices (symmetric matrices if $\H=\R$) $\A=(A_j)_{j=1}^N$ in $\H^{d\times d}$. We say $\A$ has the (generalized) {\em phase retrieval property}, or is {\em phase retrievable}, if $(\vx^*A_j\vx)_{j=1}^N$ uniquely determines $\vx\in\H^d$ up to a unimodular constant. In other words, $\vx^*A_j\vx = \vy^*A_j\vy, j=1,\ldots,N$ if and only if $\vy = c\vx$ for some $c\in\H$ and $|c|=1$. Generalized phase retrieval was studied  in  \cite{WaXu16} by Wang and Xu, and it includes the standard phase retrieval and various spinoffs as special cases. If all $A_j$ have the form $A_j = \vf_j\vf_j^*$ where $\vf_j\in\H^d$, then it is the standard phase retrieval. If all $A_j$ are orthogonal projections, namely $A_j^2=A_j$, then it becomes the fusion frame (projection) phase retrieval. Moreover, if all $A_j$ are positive semi-definite matrices satisfying $\sum_{j=1}^NA_j={\mathbb I}$ where ${\mathbb I}$ is the identity matrix, then it recasts as the POVM which is an active research topic in quantum tomography (see e.g. \cite{HMW13}).

We shall use $\Herm(\H)$ to  denote the set  of Hermitian matrices in $\H^{d\times d}$.  Just like the standard phase retrieval problem we consider the equivalence relation $\sim$ on $\H^d$: $\vx_1 \sim \vx_2$ if there is a constant $c\in \H$ with $\abs{c}=1$ such that $\vx_1=c\vx_2$. Let $\qHH :=\HH/\!\!\sim$. We shall use $\ul\vx$ to denote the equivalent class containing $\vx$. For any given $\A=(A_j)_{j=1}^N \in \Herm^N(\H)$, define the map $\MM_\A: \qHH \lra \R^N$ by
\begin{equation}  \label{eq:M-map}
      \MM_\A(\ul\vx) = (\vx^*A_1\vx, \dots, \vx^*A_N\vx)^T.
\end{equation}
Thus, the generalized phase retrieval problem asks whether a $\ul\vx\in\qHH$ is uniquely determined
by $\MM_\A(\ul\vx)$, i.e. whether $\MM_\A$ is injective on $\qHH$. We should observe that $\MM_\A$ can also be viewed as a map from $\HH$ to $\R^N$, and we shall often do this when there is no confusion.

For practical applications it is always enough to design algorithms which can recover almost all the signals \cite{FMW14,KEO16}. This gives rise to the question of almost everywhere phase retrieval. We next introduce the definition of {\em almost everywhere phase retrieval property}.

\begin{defi}  \label{defi-PRae}{\rm
    Let $\A=(A_j)_{j=1}^N \in \Herm^N(\H)$. We say $\A$ has the {\em almost everywhere phase retrieval property} or is {\em phase retrievable almost everywhere (PR-ae)}  if for almost every $\ul\vx\in\qHH$ we have $\MM_\A^{-1}\bigl(\MM_\A(\ul\vx)\bigr) = \{\ul\vx\}$ where
$\MM_\A^{-1}\bigl(\MM_\A(\ul\vx)\bigr):= \{\ul\vy\in\qHH : \MM_\A(\ul\vy)=\MM_\A(\ul\vx)\}$.
}
\end{defi}

This paper studies the following question:
\begin{itemize}
\item{\em What is the minimal $N$ for which there exists an $\A=(A_j)_{j=1}^N\in \H_d^N(\H)$ having the almost everywhere phase retrieval property? Moreover, under what conditions does $\A=(A_j)_{j=1}^N$ have the almost everywhere phase retrieval property? }
\end{itemize}
The aim of  this paper is to present a series of results addressing these questions.

\subsection{Related Work}
In the standard phase retrieval setting where $A_j = \vf_j\vf^*_j$, some answers have been provided in several studies \cite{FMW14, BCE06,FSC05, Mix14}.  For the case $\H=\R$ and $A_j=\vf_j\vf^*_j, j=1,\ldots,N$, it was shown that $N\geq d+1$  is necessary and sufficient  for there existing  $\A=(A_j)_{j=1}^N$ which is  almost everywhere  phase retrievable (see \cite{BCE06, FMW14}). Similarly, for the case $\H=\C$ and $A_j=\vf_j\vf^*_j, j=1,\ldots,N$, it is known that $2d$ generic measurements are sufficient for almost phase retrieval \cite{BCE06}.
In the context of quantum tomography, one is interested in the pure-state information complete (PSI-complete) POVM which requires POVM to determine almost all the pure state (up to a global phase) (see \cite{FSC05, F04}). In fact, PSI-complete POVM is a special case of PR-ae in which all $\{A_j\}_{ j=1}^N$ are positive semi-definite matrices satisfying $\sum_{j=1}^NA_j= {\mathbb I}$ where ${\mathbb I}$ is the identity matrix. The results in \cite{FSC05} show that there exist $N=2d$ positive semi-definite matrices $\{A_j\}_{j=1}^{2d}\subset \Herm(\C)$  satisfying $\sum_{j=1}^{2d}A_j={\mathbb I}$ which have almost everywhere phase retrieval property.  In \cite{F04} the same  result is proved with the additional requirement of  $\rank(A_j)=1$ for all $j$.

  In the generalized phase retrieval setting, these questions becomes significantly harder. For example, even in the case of fusion frames (projection) phase retrieval, the answers to these questions are far from being  known \cite{FiMi14}.

We would like to mention that the almost everywhere matrix recovery is studied in \cite{RWX}. One main result of \cite{RWX} is that $N>(p+q)r-r^2$ generic linear measurements have the almost everywhere rank-$r$ matrix recovery property in $\R^{p\times q}$ or $\C^{p\times q}$.

\subsection{Our Contribution}

In this paper we establish a general framework for the almost everywhere phase retrieval.  We prove the results for generic and random measurements  under very general settings.

For generalized phase retrieval in the real case,  we have the following theorem:

\begin{theo}  \label{theo-PRae-Generic}
   Assume that $N\geq d+1$. Then a random $\A=(A_j)_{j=1}^N$ in $\Herm^N(\R)$ chosen under any absolutely continuous distribution have the almost everywhere phase retrieval property in $\R^d$ with probability one. More generally,  let $1 \leq r_1, \dots, r_N \leq d$ and $V_j \subset \Herm(\R)$ be either the set of all rank $r_j$ symmetric matrices or the set of all rank $r_j$ orthogonal projection matrices. Then a generic $\A=(A_j)_{j=1}^N\in V_1 \times \cdots \times V_N$ has the almost everywhere phase retrieval property in $\R^d$.
\end{theo}

For generalized phase retrieval in the complex case,  we have a similar theorem:

\begin{theo}  \label{theo-ComplexGenericPRae}
   Assume that $N\geq 2d$. Then a random  $\A=(A_j)_{j=1}^N$ in $\Herm^N(\C)$ chosen under any absolutely continuous distribution has the almost everywhere phase retrieval property in $\C^d$ with probability one. More generally, let $1 \leq r_1, \dots, r_N \leq d$ and $V_j \subset \Herm(\C)$ be either the set of all rank $r_j$ Hermitian matrices or the set of all rank $r_j$ orthogonal projection matrices. Then a generic $\A=(A_j)_{j=1}^N\in V_1 \times \cdots \times V_N$ has the almost everywhere phase retrieval property in $\C^d$.
\end{theo}

Again we should point out that the above results holds for far broader classes of $V_j$ and we provide techniques for establishing PR-ae property in more general cases later in this paper.

One interesting question is whether those bounds are sharp. For almost everywhere generalized phase retrieval, there is a lower bound  $N\geq d$ for $\R^d$ and $N \geq 2d-1$ for $\C^d$ (see Corollary \ref{coro-jacobi} in Section 2). One naturally asks whether PR-ae property can be attained for $N =d$ and $N=2d-1$ in the real and complex cases, respectively. In the real case, under the standard phase retrieval assumption where all measurement matrices have rank one, this is impossible with $N \geq d+1$ being sharp. However, it is possible if the rank one property is removed . Particularly, we have the following results:


\begin{theo}\label{th:thdR}
 There exist $d$ matrices $A_1,\ldots,A_d\in \Herm(\R)$ such that $\A=(A_j)_{j=1}^d$ has the almost everywhere phase retrieval property in $\R^d$.
\end{theo}

\begin{theo}\label{th:thd}
 There exist $2d-1$ matrices $A_1,\ldots,A_{2d-1}\in \Herm(\C)$ such that $\A=(A_j)_{j=1}^{2d-1}$ has the almost everywhere phase retrieval property in $\C^d$.
\end{theo}

We would like to mention that it is possible to prove
Theorem \ref{theo-PRae-Generic} and Theorem \ref{theo-ComplexGenericPRae} using the results from \cite{RWX}. In this paper, we present a novel method for proving them. We believe that
the method developed in this paper is independent interesting and it is useful and powerful for the phase retrieval.
 For example, motivated by the methods developed in the proof of
Theorem \ref{theo-PRae-Generic} and Theorem \ref{theo-ComplexGenericPRae}, we can obtain that
Theorem \ref{th:thdR} and Theorem \ref{th:thd}. To our knowledge, Theorem \ref{th:thdR} and Theorem \ref{th:thd} are not easy to be derived from the results in \cite{RWX}.

The paper is organized as follows. In Section 2, after introducing  some notations, some of which have been used in previous studies, we also present some preliminary results on generalized almost everywhere phase retrieval, including the necessary lower bounds $N\geq d$  and $N \geq 2d-1$ for PR-ae property in $\R^d$ and $\C^d$, respectively.

In Section 3 we explore the links between phase retrieval  and the classic algebraic geometry. We recall some of the background results on the dimension
  of intersections of varieties, from which we tie the almost everywhere matrix recovery and phase
   retrieval properties to the dimension of certain varieties. These results are then used in
   Section 4 to prove the results listed above. In Section 5, we return to the standard phase retrieval with
   presenting  some additional results under this setting, which may be independent interesting.
   Particularly, we present the sufficient and necessary condition for
    $\{\vf_1,\ldots,\vf_N\}\subset \R^d$ having PR-ae property in $\R^d$. 
\section{Preliminary Results on Almost Everywhere Phase Retrieval}
\setcounter{equation}{0}

In this section, we  establish some preliminary results on almost everywhere phase retrieval. These results will play an important role for the further study of this topic.
We begin from introducing a few results and notations.
\subsection{Previous results and notations }
For any $c\in\C$ let $\Re(c)$ and $\Im(c)$ denote the real and imaginary part of $c$, respectively. A useful formula is that for a Hermitian $A\in \Herm(\H)$ and any $\vx,\vy\in\H^d$ we must have
\begin{equation} \label{2.1}
       \vx^* A\vx - \vy^*A\vy = 4\Re (\vv^*A\vu)
\end{equation}
where $\vv=\frac{1}{2}(\vx+\vy)$  and $\vu=\frac{1}{2}(\vx-\vy)$. This is straightforward to check. In the real case $\H=\R$ it means that $\vx^* A\vx - \vy^*A\vy = 4\vv^*A\vu$. In \cite{WaXu16} a series of equivalent formulations for the generalized phase retrieval have been stated. For the real case $\H=\R$ and $\A=(A_j)_{j=1}^N \in \Herm^N(\H)$, the following conclusions  are equivalent (see \cite{WaXu16}):
\begin{itemize}
   \item[\rm (1)] $\A$ has the phase retrieval property.
   \item[\rm (2)] There exist no nonzero $\vv,\vu\in\R^d$ such that $\vv^* A_j \vu=0$ for all $1\leq j\leq N$.
   \item[\rm (3)] For any nonzero $\vu\in\R^d$ we have $\Span\{A_j \vu\}_{j=1}^N=\R^d$.
   \item[\rm (4)] The Jacobian matrix of $\MM_\A$ has rank $d$ everywhere on $\R^d\setminus\{0\}$.
\end{itemize}
For the complex case $\H=\C$ and $\A=(A_j)_{j=1}^N \in \Herm^N(\H)$, the following are equivalent (see \cite{WaXu16}):
 \begin{itemize}
   \item[\rm (1)] $\A$ has the phase retrieval property.
   \item[\rm (2)] There exist no $\vv,\vu\neq 0$ in $\C^d$ with $\vu \neq ic\vv$ for any $c\in\R$ such that $\Re(\vv^* A_j \vu)=0$ for all $1\leq j\leq N$.
   \item[\rm (3)] The (real) Jacobian matrix of $\MM_\A$ has (real) rank $2d-1$ everywhere on $\C^d \setminus\{0\}$.
\end{itemize}
It is also shown in \cite{WaXu16} that the set of phase retrievable $\A$ in $\Herm^N(\H)$ is an open set, so it is stable under small perturbations.

Next, we introduce some notations which examine the set of points in $\qHH$ at which $\MM_\A$ is injective.
We shall treat $\qHH$ as a {\em real} manifold, which has dimension
$d$ if $\H=\R$ and dimension $2d-1$ if $\H=\C$.
When $\H=\R$ the Jacobian of $\MM_\A$ at $\vx$ is exactly
\[
J_\A(\vx)=2[A_1\vx, A_2\vx, \dots, A_N\vx].
\]
For the case where $\H=\C$, write $A_j = B_j+iC_j$ where $B_j, C_j$ are real. Then $B_j^\top =B_j$ and $C_j^\top =-C_j$. Let
\begin{equation}  \label{2.2}
     F_j = \begin{bmatrix}
         B_j & -C_j\\C_j & B_j
     \end{bmatrix}.
\end{equation}
For any $\vx = \Re(\vx) + i\Im(\vx)\in\C^d$ we set $\vu^\top : = [\Re(\vx)^\top , \Im(\vx)^\top ]$. Thus, the real Jacobian of $\MM_\A(\vx)$ is precisely
\begin{equation}\label{eq:cmj}
     J_\A(\vx)= 2[F_1\vu, F_2\vu, \dots, F_N\vu].
\end{equation}
 A point $\ul\vx\in\qHH$ is called
a {\em regular point} if the real Jacobian of $\MM_\A$ at
$\ul\vx$ has full rank, i.e. it has rank $d$ if $\H=\R$ and rank $2d-1$ if $\H=\C$.
Otherwise $\ul\vx$ is called a {\em degenerate point}. It is well known that
the set of all degenerate points for $\MM_\A$ is a closed set in $\qHH$.

\subsection{Preliminary Results }

\begin{theo}  \label{theo-PRaeClosure}
   For  $\A= (A_j)_{j=1}^N\in \Herm^N(\H)$, let $\Omega_\A$ be the set of points in $\ul{\H^d}$ at which
   $\MM_\A$ is not injective. Let $Z_\A$ be the set of degenerative points for $\MM_\A$
   in $\ul{\H^d}$. Then the following hold:
   \begin{itemize}
   \item[\rm (A)] $Z_\A$ is either $\ul{\H^d}$ or a null set. Furthermore if $Z_\A =\ul{\H^d}$ then $\ul{\H^d}\setminus \Omega_\A$ is a null set and hence $\A$ is not PR-ae.
   \item[\rm (B)]  Let $\overline\Omega_\A$ be the closure of $\Omega_\A$. We have $\Omega_\A\cup Z_\A \subseteq \overline\Omega_\A$. If all $A_j$ are positive semidefinite then $\Omega_\A\cup Z_\A = \overline\Omega_\A$
   \end{itemize}
\end{theo}
\Proof (A)~~ We shall prove the results for the case $\H=\C$. The real case is virtually
identical. Firstly, using the standard technique we identify the set of element $\ul\vx\in\ul{\H^d}$ with $x_1 \neq 0$ as $V_1:=\R^+\times \R^{2d-2}$.
 Now $V_1$ is almost all $\ul{\H^d}$ and its closure is $\ul{\H^d}$.
 Restricted to $V_1$ the Jacobian matrix of
$\MM_\A$ consists of entries that are linear functions (see (\ref{eq:cmj})). A point is degenerate if
and only if all $(2d-1) \times (2d-1)$ submatrices of the Jacobian matrix at that
point have determinants 0. Note that each determinant is a polynomial. Thus, the set of
degenerate points in $V_1$ is the intersection of real algebraic varieties in
$\R^{2d-1}$ restricted to $V_1=\R^+\times \R^{2d-2}$. It follows that
the set is either all $V_1$, or a null set with local dimension less than
$2d-1$. Hence,
$Z_\A$ is either $\overline{V_1}=\ul{\H^d}$ or a null set.

Assume that $Z_\A = \ul{\H^d}$. Then the Jacobian of $\MM_\A$ has maximal rank strictly less than $2d-1$ at any point.
 Let $\Gamma$ be the
set of points in $\ul{\H^d}$ at which the Jacobian of $\MM_\A$ has the maximal rank, say $r$.
Then $\Gamma$ is an open set. Furthermore, the complement set of $\Gamma$ are
precisely the points at which all $(r-1)\times (r-1)$ submatrices have
zero determinant. Thus, $\Gamma^c$ is a null set by the same argument as before.
The Rank Theorem now implies that the map $\MM_\A$ is not injective at a
neighborhood of any point in $\Gamma$. Thus, $\ul{\H^d}\setminus \Omega_\A$ is a null set,
i.e. almost all points in $\ul{\H^d}$ are not injective for $\MM_\A$. Thus, $\A$ is not PR-ae.

\vspace{1ex}

(B)~~ Again we only need to consider the case $\H=\C$. The real case is virtually identical.
We first prove that $Z_\A\subseteq \overline\Omega_\A$. For any $\ul\vx\in Z_\A$ the rank of the (real) Jacobian $J_\A$ of $\MM_\A$ is at most $2d-2$ at $\vx \in\H^d$. Thus, there exists a $\vv\in \H^d$ such that $\vv \neq ic \vx$ and $\Re(\vv^*J_\A(\vx))=\Re(\vv^*A_j\vx )=0$ where $c\in \R$. Let $t_k=1/k$ and set $\vx_k =\vx+t_k\vv$, $\vy_k = \vx-t_k\vv$. Then $\vx_k \neq \vy_k$ in $\qHH$. Combining $\Re(\vv^*A_j\vx )=0$  and (\ref{2.1}) we obtain  $\MM_\A(\ul\vx_k)=\MM_\A(\ul\vy_k)$, which implies that $\ul\vx_k, \ul\vy_k \in \Omega_\A$. Clearly $\lim_k \ul\vx_k = \ul\vx$, and thus, $\ul\vx \in\overline\Omega_\A$.

Now assume that all $A_j$ are positive semi-definite matrices. We first prove $\overline\Omega_\A \subseteq \Omega_\A\cup Z_\A$.
For each $\delta>0$ define $E_\delta$ to be the set of all $\ul\vx\in\ul{\H^d}$
such that there exists a $\ul\vy \in\ul{\H^d} $ with $\MM_\A(\ul\vy)=\MM_\A(\ul\vx)$ and
\[
\uD(\vx,\vy):=\min_{\alpha\in \H,\abs{\alpha}=1}\|\vx-\alpha\vy\|\,\,\geq\,\, \delta.
 \]
 The definition of $E_\delta$ implies  that it  is a closed set in
$\ul{\H^d}$. Choose a positive sequence $\delta_k\downarrow 0$. Then
$\Omega_\A = \bigcup_{k} E_{\delta_k}$.

Let $\ul\vx_k$ be a sequence in $\Omega_\A$ with $\ul\vx_k \lra \ul\vx$. We need to show
$\ul\vx\in \Omega_\A\cup Z_\A$. If $\ul\vx=0$ we arrive at the conclusion since $0\in Z_\A$.
The conclusion also holds  if $\ul\vx\in\Omega_\A$.  It remains to consider the case where $\ul\vx\not\in \{0\}\cup \Omega_\A$. Because each $E_{\delta_k}$ is closed, without loss of generality,  we  assume that $\ul\vx_k \in E_{\delta_k}$. Let $\ul\vy_k \in\MM_\A^{-1}(\MM_\A(\ul\vx_k))$
such that $\uD(\ul\vx_k, \ul\vy_k) \geq \delta_k$.
We first consider the case where $\{ \ul\vy_k\}_{k\in \Z}$ is a bounded sequence. Then there exist 
an accumulation point $\ul\vy$ and a subsequence where we still denote it by $\{ \ul\vy_k\}_{k\in \Z}$ such that  $\ul\vy_k \lra \ul\vy$. Clearly $\ul\vy\in\MM_\A^{-1}(\MM_\A(\ul\vx))$. Recall the assumption $\ul\vx\not\in \{0\}\cup \Omega_\A$, which implies that $\ul\vy = \ul\vx$. Thus,  any
small neighborhood of $\ul\vx$ contains $\ul\vx_k$ and $\ul\vy_k$ for sufficiently
large $k$, which means $\MM_\A$ is not one to one locally at $\ul\vx$. However, a smooth map must be locally one to one at a regular point. It means that $\ul\vx$ is not a regular point of $\MM_\A$, i.e. $\ul\vx\in Z_\A$.

We still need to consider the case where $\{ \ul\vy_k\}_{k\in \Z}$ is unbounded.
 Set
 \[
 {\mathcal N}:=\{\eta\in \H^d: \eta^*A_j\eta=0, j=1,\ldots,N\}.
 \]
 Since $A_1,\ldots,A_N$ are positive semi-definite, ${\mathcal N}\subset \H^d$ is a linear space.
 We take $\vy_k=\vy'_k\bigoplus \eta_k$ where $\eta_k\in  {\mathcal N}$. Then $\{\vy'_k\}_{k\in \Z}$ must be a bounded sequence (otherwise, $\MM_A(\vy'_k)=\MM_\A(\vy_k)$ is a unbounded sequence).
 Note that   $\MM_A(\vy'_k)=\MM_\A(\vy_k)$. Then we can replace $\ul\vy_k$ by $\ul\vy'_k$ in the argument above and obtain that $\ul\vx\in Z_\A$.
\eproof

\vspace{2mm}

\begin{coro}  \label{coro-jacobi}
   Assume that $\A=(A_j)_{j=1}^N\in \Herm^N(\H)$ is PR-ae on $\H^d$. Then $N\geq d$ for $\H=\R$  and $N\geq 2d-1$ for $\H=\C$.
\end{coro}
\Proof
 According to the definition of degenerative  points, $Z_\A=\ul{\H^d}$ if $N<d$ for $\H=\R$ and $N<2d-1$ for $\H=\C$. According to  the (A) of Theorem \ref{theo-PRaeClosure},    $\A=(A_j)_{j=1}^N$ is not PR-ae on $\H^d$. We arrive at the conclusion.
\eproof


We  introduce the following lemma which plays an important role in this paper.
\begin{lem}  \label{prop-PositiveMeasure}
   Let $\A=(A_j)_{j=1}^N\in \Herm^N(\H)$ where $\H=\R$ or $\C$. The following are equivalent:
   \begin{itemize}
   \item[\rm (i)] $\A$ is not PR-ae.
   \item[\rm (ii)] Let ${\mathcal U}$ be the set of $(\vu,\vv)\in \H^d\times\H^d$
   with $\vu \neq ic\vv$ for any $c\in\R$ such that
   \begin{equation*}
       \Re (\vv^*A_j\vu)=0 \mhsp  \mbox{for all $1\le j\le N$.}
   \end{equation*}
   The set $E=\{\vu+\vv:~(\vu,\vv)\in {\mathcal U},\vu\neq 0,\vv\neq 0 \}$ has positive
    Lebesgue outer   measure in $\H^d$.
   \end{itemize}
\end{lem}
\Proof  To see {\rm (ii)} $\Rightarrow$ {\rm (i)}, for any $\vx=\vu+\vv\in E$ with $(\vu,\vv)\in {\mathcal U}$, set $\vy = \vu-\vv$. Then according to (\ref{2.1}), one has $\MM_\A(\vx)=\MM_\A(\vy)$. Since $(\vu,\vv)\in {\mathcal U}$, it implies that $\vx \neq \alpha\vy$ for any $\abs{\alpha}=1$. Note that $E$ has positive measure, which gives that $\MM_\A$ is not injective
 in a set with positive Lebesgue outer measure and hence $\A$ is not PR-ae. The converse {\rm (i)} $\Rightarrow$ {\rm (ii)} follows from the similar argument.
\eproof

\vspace{3mm}

\begin{remark} {\rm
For the case where $\H=\R$,  the {\rm (ii)} in Lemma \ref{prop-PositiveMeasure} is reduced to the following statement: Suppose that ${\mathcal U}$ is the set of $(\vu,\vv)\in \H^d\times\H^d$ with $\vu\neq 0, \vv\neq 0$  such that
   \begin{equation*}
       \vv^*A_j\vu=0 \mhsp  \mbox{for all $1\le j\le N$.}
   \end{equation*}
   The set $E=\{\vu+\vv:~(\vu,\vv)\in {\mathcal U}\}$ has positive
    Lebesgue  measure in $\R^d$.
}
\end{remark}

\section{ Almost Everywhere Phase Retrieval and The Dimension of Algebraic Variety}
\setcounter{equation}{0}

The phase retrieval problem has a well known formulation in terms of low rank matrices \cite{CESV12,CSV12, Xu15}. Particularly, for any $A\in\H^{d\times d}$ and $\vx,\vy\in\H^d$, it holds that $\vy^T A\vx=\tr(AQ)$ where $Q=\vx\vy^T$. This relationship transforms phase retrieval into a  recovering of a rank-one matrix. In this section, we extend this relationship further (see Theorem \ref{theo-PRae-Variety} and  Theorem \ref{theo-PRae-VarietyC}).
Before proceeding to these results, we first introduce some basic notations related to projective spaces and varieties.

 \subsection{Background from algebraic geometry}
 For any complex vector space $X$ we shall use $\PP(X)$ to denote the induced projective space, i.e. the set of all one dimensional subspaces in $X$. As usual for each $\vx\in X$ we use $[\vx]$ to denote the induced elements in $\PP(X)$. Similarly, for any subset $S\subset X$ we use $[S]$ or $\PP(S)$ to denotes its projectivization in $\PP(X)$. Throughout this paper, we say $V\subset \C^d$ is a projective variety if $V$ is the locus of a collection of homogeneous  polynomials in $\C[\vx]$.
Strictly speaking a projective variety lies in $\PP(\C^d)$ and is the projectivization of the zero locus of a collection of homogeneous polynomials. But like in \cite{WaXu16}, when there is no confusion the phrase {\em projective variety in $\C^d$} means an algebraic variety in $\C^d$ defined by homogeneous polynomials. We shall use {\em projective variety in $\PP(\C^d)$} to describe a true projective variety. Note that sometimes it is useful to consider the more general quasi-projective varieties. A set $U \subset \C^d$ is a {\em quasi-projective variety} if there exist two projective varieties $V$ and $Y$ with $Y\subset V$ such that $U=V\setminus Y$.

We shall use $V\cap \R^d$ to denote the real points of $V$. A key fact is that for a variety $V$ we have $\dim_\R(V\cap \R^d)\leq \dim(V)$ (see Section 2.1.3 in \cite{E15} and \cite{WaXu16}). This also holds for a quasi-projective variety since the proof uses only local properties of $V$ (see \cite{WaXu16}). The definitions of $\dim(V)$ and $\dim_\R(V\cap \R^d)$ are introduced in \cite{alge} and \cite{real}, respectively (see also  \cite{E15}).

In this paper we shall often focus on studying the set
\begin{equation}  \label{eq:rankr}
  {\mathcal M}_{d,r}(\H):=\Bigl\{Q\in \H^{d\times d}: {\rank}(Q)\leq r\Bigr\}, \mhsp\H=\C \mbox{~or~}\R.
\end{equation}
Note that ${\rank}(Q)\leq r$ is  equivalent to the vanishing of all $(r+1)\times (r+1)$ minors of $Q$.  Hence, ${\mathcal M}_{d,r}(\H)$ is a well-defined projective variety in
$\H^{d\times d}$ with  $\dim_\H({\mathcal M}_{d,r}(\H))=2dr-r^2$ \cite[Prop. 12.2]{alge}. More generally, for the matrix recovery problem we will consider non-square $p$ by $q$ matrices of rank $r$:
\begin{equation}  \label{eq:rankr-rectangle}
  {\mathcal M}_{p\times q,r}(\H):=\Bigl\{Q\in \H^{p\times q}: {\rank}(Q)\leq r\Bigr\}, \mhsp\H=\C \mbox{~or~}\R.
\end{equation}
Again, ${\mathcal M}_{p\times q,r}(\H)$ is a projective variety and it is well known that
$\dim_\H {\mathcal M}_{p\times q,r}(\H) = r(p+q)-r^2$.

In \cite{WaXu16} the notion of an admissible algebraic variety with respect to a family of linear functions was introduced. This concept is equally useful in this paper.

\begin{defi}[\cite{WaXu16}]   \label{defi-admissible}  {\rm
   Let $V$ be the zero locus of a finite collection of homogeneous polynomials in $\C^{M}$ with $\dim V>0$ and let $\{\ell_\alpha(\vx): \alpha\in I\}$ be a family of (homogeneous) linear functions. We say $V$ is {\em admissible} with respect to $\{\ell_\alpha(\vx)\}$ if  $\dim (V\cap \{\ell_\alpha(\vx)=0\}) <\dim V$ for all $\alpha\in I$.
}
\end{defi}

It is well known in algebraic geometry that if $V$ is irreducible in $\C^M$ then $\dim (V\cap Y) = \dim(V)-1$ for any hyperplane $Y$ that does not contain $V$. Thus, the above admissible condition is equivalent to the property that no irreducible component of $V$ of dimension $\dim (V)$ is contained in any hyperplane $\ell_\alpha(\vx)=0$. In general without the irreducibility condition, admissibility is equivalent to that for a generic point $\vx\in V$, any small neighborhood $U$ of $\vx$ has the property that $U\cap V$ is not completely contained in any hyperplane $\ell_\alpha(\vx)=0$.

 Many projective varieties have the required admissibility property. We just list a few of them below:

\begin{prop}( \cite[Proposition 4.1]{RWX})  \label{prop:adm}
     Let $V$ be one of the following projective varieties in $\C^{q\times p}$. Then $V$  is admissible with respect to the maps $\{ \phi_Q(\cdot )=\tr(\cdot Q):~Q\in {\mathcal M}_{p\times q,r}(\C)\}$, where $1 \leq r \leq \frac{1}{2} \min (p,q)$:
\begin{itemize}
\item[\rm (A)]~~$V =  {\mathcal M}_{q\times p,s}(\C)$, where $1 \leq s \leq \min(p,q)$.
\item[\rm (B)]~~$q \geq p$ and $V$ is the set of all scalar multiples of matrices $P$ whose rows are complex orthonormal in the sense that any two rows $\vx,\vy$ of $P$ have $\vx\vy^\top =\delta(\vx-\vy)$.
\item[\rm (C)]~~$q \leq p$ and $V$ is the set of  all scalar multiples of matrices whose columns are complex orthonormal in the sense that any two columns $\vx,\vy$ of $P$ have $\vx^\top\vy =\delta(\vx-\vy)$.
\item[\rm (D)]~~$q = p=d$ and $V$ is the set of all  all scalar multiples of $d\times d$ rank $s$ complex orthogonal projection matrices in the sense that $P=P^*$ and $P^2=P$.
\end{itemize}
\end{prop}

  The following proposition is  from
\cite{RWX}:

\begin{theo} ( \cite[Theorem 2.1]{RWX}) \label{theo-admissible}
For $j=1, \dots, N$ let $L_j:\C^{n} \times \C^{m}\rightarrow \C$ be bilinear functions  and $V_j$ be projective varieties in $\C^n$.
Set $V := V_1\times \dots \times V_N \subseteq (\C^{n})^N$. Let $W, Y\subset \C^{m}$ be a projective variety in $\C^m$, $W \setminus Y$ be a quasi-projective variety.
For each fixed $j$, assume that $V_j$ is admissible with respect to the linear functions $\{f^\vw (\cdot)=L_j(\cdot,\vw):~\vw\in W\setminus Y\}$.
\begin{itemize}
\item[(1)]~~Assume that $N \geq \dim W$. There exists an algebraic subvariety $Z \subseteq V$ with $\dim (Z) < \dim (V)$ such that for any $\vx =(\vv_j)_{j=1}^N \in V\setminus Z$, the subvariety  $X_{\vx}$ given by
$$
    X_{\vx} := \Bigl\{\vw\in W \setminus Y:~L_j(\vv_j,\vw) = 0 \mbox{~for all $1\leq j \leq N$}\Bigr\}
$$
is the empty set.
\item[(2)]~~Assume that $N < \dim W$. There exists an algebraic subvariety $Z\subset V$  with $\dim Z < \dim V$ such that for any $\vx =(\vv_j)_{j=1}^N \in V\setminus Z$,  the subvariety  $X_{\vx}$ given by
$$
    X_{\vx} := \Bigl\{\vw\in W \setminus Y:~L_j(\vv_j,\vw) = 0 \mbox{~for all $1\leq j \leq N$}\Bigr\}
$$
has $\dim X_{\vx} = \dim W-N$.
\end{itemize}
\end{theo}


\medskip

\subsection{Almost Everywhere Phase Retrieval: Real Case}
In this subsection, we consider the almost everywhere phase retrieval for the case where $\H=\R$.

\begin{theo}  \label{theo-PRae-Variety}
Assume that $\A=(A_j)_{j=1}^N\in \Herm^N(\R)$. Let $X_\A \subset \C^{d\times d}$ be given by
\begin{equation}  \label{3.1}
    X_\A := \Bigl\{Q\in\C^{d\times d}:~\rank(Q)\leq 1 \mbox{~\rm{and}~} \tr(A_jQ)=0 \mbox{~\rm{for all }$1\leq j \leq N$}\Bigr\}.
\end{equation}
If the (complex) variety $X_\A$  has dimension $\dim(X_\A) \leq d-2$, then $\A$ has the PR-ae property in $\R^d$.
\end{theo}
\Proof  Because $X_\A$ is the zero locus of some homogeneous polynomials, we can view it naturally as a projective variety with $\dim(\PP(X_\A)) = \dim(X_\A)-1 \leq d-3$.  Consider the map $\Phi: \PP(\C^{d})\times \PP(\C^{d}) \lra \PP(\C^{d\times d})$ given by
$$
      \Phi ([\vx], [\vy]) = [\vx\vy^T].
$$
It is easy to see that $\Phi$ is injective, and furthermore it is a one-to-one mapping of $\PP(\C^{d})\times \PP(\C^{d})$ onto the set of rank one matrices in $\PP(\C^{d\times d})$. Note that the projective variety $[X_\A] =\PP(X_\A) \subset \PP(\C^{d\times d})$ has dimension $\dim([X_\A])=\dim(X_\A)-1 \leq d-3$. Hence, the dimension of the projective variety $\Phi^{-1}([X_\A])$ has dimension at most $d-3$. Let
$$
     Y_\A:= \Bigl\{(\vx,\vy)\in \C^{d} \times \C^d:~\Phi([\vx], [\vy])\in [X_\A]\Bigr\}.
$$
It follows that $\dim (Y_\A) =\dim([X_\A])+2 \leq d-1$.

Noting $\vx^TA_j\vy={\rm tr}(A_j\vy\vx^T)$, we have $\vx^TA_j\vy=0$ provided $
(\vx, \vy)\in Y_\A\cap\R^d\times \R^d$.
According to Lemma \ref{prop-PositiveMeasure}, to this end, it is enough to show that $\{\vx+\vy: \vx, \vy\in Y_\A\cap\R^d\times \R^d\}$ has zero Lebesgue measure.
Note that the real slice of a complex algebraic variety is a real algebraic variety whose real dimension is no more than the dimension of the complex variety (see \cite[Lemma 3.1]{WaXu16}). Thus, the real slice of $Y_\A$ has real dimension $\dim_\R(Y_\A\cap\R^d\times \R^d) \leq \dim Y_\A \leq d-1$. Let $\tau: \R^d\times\R^d \lra \R^d$ be given by $\tau(\vx,\vy) = \vx+\vy$. Then $\dim_\R (\tau(Y_\A\cap\R^d\times \R^d)) \leq d-1$ and thus,
\[
\tau(Y_\A\cap\R^d\times \R^d)= \{\vx+\vy: \vx, \vy\in Y_\A\cap\R^d\times \R^d\}
\]
 has zero Lebesgue measure. By Lemma \ref{prop-PositiveMeasure}, $\A$ must have the PR-ae property.
\eproof

\subsection{Almost Everywhere Phase Retrieval: Complex Case}
We now turn our attention to the complex setting $\H=\C$. First we establish a couple of auxiliary results.

\begin{lem}   \label{lem:matrix-indep}
   Assume that  $\vx,\vy\in\C^d$ are linearly independent. Then $\vx\vx^*, \vy\vy^*, \vx\vy^*, \vy\vx^*$ are linearly independent in $\C^{d\times d}$.
\end{lem}
\Proof  Since $\vx,\vy$ are linearly independent, there exists a $\vv\in\C^d$ such that $\vx^*\vv=1$ and $\vy^*\vv=0$. Assume that
\begin{equation}\label{eq:ind}
     c_1 \vx\vx^* +c_2 \vy\vy^*+c_3 \vx\vy^*+c_4 \vy\vx^* = 0
\end{equation}
where $c_1, \dots, c_4\in\C$. Then multiplying $\vv$ on the both sides of (\ref{eq:ind})  we obtain
$$
      c_1 \vx\vx^*\vv +c_2 \vy\vy^*\vv+c_3 \vx\vy^*\vv+c_4 \vy\vx^*\vv =
       c_1 \vx++c_4 \vy = 0.
$$
Hence, $c_1=c_4 = 0$. Similarly, there exists a $\vu\in\C^d$ such that $\vy^*\vu=1$ and $\vx^*\vu=0$. Multiplying $\vu$ on the both sides yields
$$
      c_1 \vx\vx^*\vu +c_2 \vy\vy^*\vu+c_3 \vx\vy^*\vu+c_4 \vy\vx^*\vu =c_2\vy +c_3\vx = 0.
 $$
Thus, $c_2=c_3=0$. The lemma is proved.
\eproof

\begin{lem}   \label{lem:rank2-matrix}
   Let  $\vx,\vy\in\C^d$ be linearly independent. Assume that $\vz\vz^*-\vw\vw^*=\lambda\vx\vx^*-\mu\vy\vy^*$ where $\vz,\vw\in\C^d$ and $\lambda,\mu\geq 0$  are not all $0$. Then $\vz,\vw\in \Span\{\vx,\vy\}$. Furthermore, set
   $$
       \vz = a\vx+b\vy,~\vw= s\vx+t\vy, \shsp \mbox{where}  \shsp a,b,s,t\in\C.
   $$
Then $\vz\vz^*-\vw\vw^*=\lambda\vx\vx^*-\mu\vy\vy^*$ if and only if there exist $\omega_1,\omega_2, \omega_3\in\C$ with $|\omega_j|=1$ and $0 \leq \beta<1$ such that
\begin{equation}    \label{rank2-formula}
     \vz = \frac{\omega_1\lambda}{\sqrt{1-\beta^2}}\vx + \frac{\omega_2\beta\mu}{\sqrt{1-\beta^2}}\vy, ~~
     \vw =  \frac{\omega_3\beta\lambda}{\sqrt{1-\beta^2}}\vx
           +\frac{\bar\omega_1\omega_2\omega_3\mu}{\sqrt{1-\beta^2}}\vy.
\end{equation}
\end{lem}
\Proof
We first consider the case where $\vz,\vw$ are linearly dependent. Then $\vz\vz^*-\vw\vw^*$ has rank at most 1. Note that if $\lambda \neq 0,\mu\neq 0$ then $\lambda\vx\vx^*-\mu\vy\vy^*$ has rank 2. This is a contradiction. Thus, in this case we must have $\lambda=0$ or $\mu=0$, say $\lambda\neq 0$ but $\mu=0$. It follows that $\vz,\vw$ must be colinear with $\vx$. Hence, $\vz,\vw\in \Span\{\vx,\vy\}$.

Now assume that $\vz,\vw$ are linearly independent.  Set
 \[
 H_{\vx,\vy}:=\{\vv\in \C^d: \vx^*\vv=0, \vy^*\vv=0\}.
 \]
 A simple observation is that $H_{\vx,\vy}$ is a linear space with $\dim(H_{\vx,\vy})=d-2$.
 The definition of $H_{\vx,\vy}$ implies that
 \[
 H_{\vx,\vy}^\bot\,\,=\,\,\Span\{\vx,\vy\}.
 \]
 For any $\vv\in H_{\vx,\vy}$  we clearly have
$$
    \vz\vz^*\vv-\vw\vw^*\vv=\lambda\vx\vx^*\vv-\mu\vy\vy^*\vv = 0.
$$
Since $\vz,\vw$ are linearly independent, we must have $\vz^*\vv=\vw^*\vv=0$. Hence, $\vv$ must be orthogonal to both $\vz,\vw$. It implies that $\vz,\vw\in  H_{\vx,\vy}^\bot= \Span\{\vx,\vy\}$.

Finally, for $\vz = a\vx+b\vy$ and $\vw= s\vx+t\vy$ we have
\begin{eqnarray*}
      \vz\vz^*-\vw\vw^* &=& (|a|^2-|s|^2)\vx\vx^* -(|t|^2-|b|^2)\vy\vy^*
             +(a\bar b-s\bar t)\vx\vy^* +(b\bar a-t\bar s)\vy\vx^*.
\end{eqnarray*}
It follows from Lemma \ref{lem:matrix-indep} that $\lambda^2\vx\vx^*-\mu^2\vy\vy^* = \vz\vz^*-\vw\vw^*$ if and only if $|a|^2-|s|^2= \lambda^2$, $|t|^2-|b|^2=\mu^2$ and $a\bar b-s\bar t=0$. Since one of $\lambda$ or $\mu$ is nonzero, say $\lambda>0$. Set $\beta = \abs{s}/\abs{a}$, then $\beta\in [0,1)$. Moreover, since $\bar b/\bar t = s/a$, it implies that $\abs{b} = \beta \abs{t}$. Hence, $\lambda^2 = (1-\beta^2)\abs{a}^2$ and $\mu^2 = (1-\beta^2)\abs{t}^2$ which gives that $\abs{a}=\frac{\lambda}{\sqrt{1-\beta^2}}, \abs{t}=\frac{\mu}{\sqrt{1-\beta^2}}$. Noting that $\abs{b}=\beta\abs{t}, s= a\bar b/\bar t $ we obtain
$$
    a = \frac{\omega_1\lambda}{\sqrt{1-\beta^2}}, ~b=\frac{\omega_2\beta\mu}{\sqrt{1-\beta^2}},
    ~s = \frac{\omega_3\beta\lambda}{\sqrt{1-\beta^2}},~ t = \frac{\omega_4\mu}{\sqrt{1-\beta^2}}
$$
for some $\omega_j\in \C$ with $|\omega_j |=1$, $1 \leq j \leq 4$. The $a\bar b-s\bar t=0$ implies that $\omega_4=\bar\omega_1\omega_2\omega_3$. The lemma is proved.
\eproof

\begin{theo}  \label{theo-PRae-VarietyC}
Assume that $\A=(A_j)_{j=1}^N\in \Herm^N(\C)$. Let $X_\A \subset \C^{d\times d}$ be given by
\begin{equation}  \label{3.1}
    X_\A := \Bigl\{Q\in {\mathcal M}_{d,2}(\C):~ \tr(A_jQ)=0 \mbox{~for all $1\leq j \leq N$}\Bigr\}.
\end{equation}
If the (complex) variety has dimension $\dim(X_\A) \leq 2d-4$, then $\A$ has the PR-ae property in $\C^d$.
\end{theo}
\Proof  Let
 \[
F=\{(\vx,\vy)\in \C^d\times \C^d:~\ul\vx\neq \ul\vy, ~\vx^*A_j\vx=\vy^*A_j\vy  \mbox{~for all~} 1\leq j \leq N\}.
 \]
 Then $\A$ has the PR-ae property if and only if $\pi_1(F)$ has zero Lebesgue outer measure, where $\pi_1((\vx,\vy))=\vx$.  Let
$$
        Y_\A = \Bigl\{ \vx^*\vx-\vy^*\vy:~(\vx,\vy) \in F\Bigr\}.
$$
Then we have $Y_\A \subset X_\A \cap \Herm(\C)$ since $\vx^*A_j\vx-\vy^*A_j\vy = \tr(A_j(\vx\vx^*-\vy\vy^*))$.

We identify $\C^M$ with $\R^M\times \R^M$ and $\C^{d\times d}$ with $\R^{d\times d} \times \R^{d\times d}$. By doing so we next show that $ X_\A \cap \Herm(\C)$ has real dimension at most $2d-4$ provided $\dim(X_\A) \leq 2d-4$. Consider the linear map $\varphi$ on $\C^{d\times d}$ given by
$$
     \varphi(A) = \frac{1}{2}(A+A^T) + \frac{i}{2}(A-A^T).
$$
It is easy to check that $\varphi$ is an isomorphism on $\C^{d\times d}$ with inverse map $\varphi^{-1}(B) =\frac{1}{2} (B+B^T)-\frac{i}{2}(B-B^T)$. Hence, $\varphi^{-1}(X_\A)$ is a variety with the same dimension as $X_\A$. It gives that the real slice $\varphi^{-1}(X_\A)\cap \R^{d\times d}$ must have real dimension no more than $2d-4$, i.e., $\dim(\varphi^{-1}(X_\A)) \leq 2d-4$. Therefore,
$$
    \dim_\R(\varphi( \varphi^{-1}(X_\A)\cap \R^{d\times d})) = \dim_\R(X_\A \cap \varphi(\R^{d\times d})) \leq 2d-4.
$$
Note that $\varphi(\R^{d\times d})$ is precisely $\Herm(\C)$. Hence,
$$
     \dim_\R(X_\A \cap \Herm(\C))\,\, \leq\,\, 2d-4.
$$

Recall the definition of the equivalence relation $\sim$ on $\C^d$ where $\vx \sim \vy$ if and only if $\vx=c\vy$ for some $c\in\C$ and $|c|=1$. For each $\vx\in\C^d$, the equivalent class containing $\vx$ is denoted by $\ul \vx$. Consider the subset $\Lambda$ of $\ul\C^d\times \ul\C^d\times\R_+^2$ given by
$$
      \Lambda = \Bigl\{(\ul\vx, \ul\vy, \lambda, \mu):~\|\vx\|=\|\vy\|=1,  \vx^*\vy=0, \lambda\geq 0, \mu\geq 0\Bigr\}.
$$
Define the map $\Psi$ on $\Lambda$ by
$$
     \Psi (\ul\vx, \ul\vy, \lambda, \mu) =  \lambda^2\vx\vx^* - \mu^2\vy\vy^*.
$$
We claim that $\Psi$ is injective. To see this we first observe that $\vx\vx^* = \vx_1\vx_1^*$ if and only if $\vx \sim\vx_1$. Assume that $\Psi(\ul\vx, \ul\vy, \lambda, \mu)=\Psi(\ul\vz, \ul\vw, \hat\lambda, \hat\mu)$. Then
$$
   \lambda^2\vx\vx^* - \mu^2\vy\vy^* = \hat\lambda^2\vz\vz^* - \hat\mu^2\vw\vw^*.
$$
Note that the eigenvectors of $\lambda^2\vx\vx^* - \mu^2\vy\vy^*$ are $\vx, \vy$ with corresponding to eigenvalues $\lambda^2, -\mu^2$, respectively. We must have $\vz=c_1\vx$, $\vw=c_2\vy$ with $|c_1|=|c_2|=1$ and $\lambda^2=\hat\lambda^2$, $\mu^2=\hat\mu^2$. Hence, $\lambda=\hat\lambda$, $\mu=\hat\mu$ and $\ul\vx=\ul\vz$, $\ul\vy = \ul\vw$.

Next we claim that the range of $\Psi$ contains $Y_\A$, i.e., $Y_\A\subset \Psi(\Lambda)$. Indeed,  for every $B\in Y_\A$, it can be decomposed into $B = \sum_{j=1}^r \lambda_j \vv_j\vv_j^*$ where $r$ is the rank of $B$ and $\{\vv_j\}_{j=1}^r$ are orthonormal. Since the rank of $B$ is at most 2 and the top two eigenvalues of $B$ cannot be both positive or negative, we can write $B$ in the form of $B= \lambda_1\vv_1\vv_1^* + \lambda_2\vv_2\vv_2^*$ with $\lambda_1 \geq 0$ and $\lambda_2 \leq 0$. So $B=\Psi(\ul\vv_1,\ul\vv_2,\lambda_1,-\lambda_2)$ which implies that $Y_\A\subset \Psi(\Lambda)$.

Now define $\Gamma = \{(\omega_1,\omega_2,\beta)\in \C^2\times \R:~\abs{\omega_1}=\abs{\omega_2}=1\}$. Then $\dim_\R\Gamma =3$. Hence,
$$
      \dim_\R(\Psi^{-1}(Y_\A) \times\Gamma) \leq 2d-4+3 = 2d-1.
$$
We shall identify each element $\ul\vx \in \ul\C^d$ with the unique element $\hat\vx$ in the equivalent class $\ul\vx$ with the property that the first nonzero entry of $\hat\vx$ is real and positive. Define the map $\pi:  \Psi^{-1}(Y_\A) \times\Gamma \lra \C^d$ by
$$
      \pi \Bigl((\ul\vx, \ul\vy, \lambda, \mu), (\omega_1,\omega_2,\beta)\Bigr)
         = \frac{\lambda\omega_1}{\sqrt{1-\beta^2}}\hat\vx+\frac{\omega_2\beta \mu}{\sqrt{1-\beta^2}}\hat \vy.
$$
By Lemma \ref{lem:rank2-matrix} the set $ \pi(\Psi^{-1}(Y_\A) \times\Gamma)= \pi_1(F)$. However, the real dimension of
$ \pi(\Psi^{-1}(Y_\A) \times\Gamma)$ is bounded from above by the real dimension of $\Psi^{-1}(Y_\A) \times\Gamma$, which is at most $2d-1$. Hence, $\pi_1(F)$ cannot have positive Lebesgue measure in $\C^d$. This means $\A$ has the PR-ae property.
\eproof

\section{Proofs of Main Results}
\setcounter{equation}{0}

In this section we apply the results from the previous sections to obtain more concrete results for phase retrieval. Particularly, we present the proofs of Theorems which stated in Section 1.

\medskip
We first consider the almost everywhere phase retrieval in the real case. Theorem \ref{theo-PRae-Generic} shows that $N\geq d+1$ generic matrices have almost everywhere phase retrieval property in $\R^d$. We now present the proof of it.

\medskip

\noindent
{\bf Proof of Theorem \ref{theo-PRae-Generic}:}~~First we consider the case of rank $r_j$ symmetric matrices. We slightly abuse the notation by extending $V_j$ to be the set of symmetric matrices in $\C^{d\times d}$ with rank no more than  $r_j$.  For any $A,Q \in \C^{d\times d}$ define $L(A,Q):=\tr(AQ)$. Note that $\dim_\R((V_j)_\R)=\dim(V_j)$.
Thus, according to Theorem 4.1 in \cite{WaXu16}, $V_j$ is admissible with respect to the linear functions $\{f^Q(A)= L(A,Q) :~Q\in {\mathcal M}_{d,1}(\C)\}$.
This implies, through Theorem \ref{theo-admissible}, that a generic real $\A\in V_1 \times \cdots \times V_N$ has dimension $\dim(X_\A) = 2d-1-N \leq d-2$, where
$$
      X_{\A}  := \Bigl\{Q\in {\mathcal M}_{d,1}(\C):~L(A_j,Q) = 0 \mbox{~for all $1\leq j \leq N$}\Bigr\}.
$$
Thus, a generic real $\A$ has the almost everywhere PR property by Theorem \ref{theo-PRae-Variety}.

%

For the case of orthogonal projection matrices, let $U_j$ be the set of all scalar multiples of complex orthogonal projection matrices  in the sense $\rank(P)=r_j$, $P=P^\top$ and $P^2=P$. Applying the exact same arguments as before, we can prove the theorem for real orthogonal projections.

Finally, the random case is a direct corollary of the first case with all $r_j = d$.
\eproof

To guarantee   $\A=(A_j)_{j=1}^N$ having PR-ae property, Theorem \ref{theo-PRae-Generic} requires that $N\geq d+1$. An interesting question is whether it is possible to lower $N$ to $d$. Theorem \ref{th:thdR} shows it is possible. We next prove that theorem.

\medskip
\noindent
{\bf Proof of Theorem \ref{th:thdR}:}~~
Suppose that the $(j,k)$ elements of $A_t\in \Herm(\R)$ are $1$ if $j+k=t+1$ and other elements are $0$, i.e.,
\[
     A_1 = \begin{bmatrix}
         1 & 0 &\cdots  &0 \\
         0 & 0 &\cdots &0  \\
         \vdots & \vdots &\cdots &\vdots  \\
         0 & 0 &\cdots &0  \\
     \end{bmatrix},\,\,
       A_2 = \begin{bmatrix}
         0 & 1 &\cdots  &0 \\
         1 & 0 &\cdots &0  \\
         \vdots & \vdots &\cdots &\vdots  \\
         0 & 0 &\cdots &0  \\
     \end{bmatrix},\ldots, A_d = \begin{bmatrix}
         0 & 0 &\cdots  &1 \\
         0 & 0 &\cdots &0  \\
         \vdots & \vdots &\cdots &\vdots  \\
         1 & 0 &\cdots &0  \\
     \end{bmatrix}.
\]

Suppose that $\vv=(v_1,\ldots,v_d),\vu=(u_1,\ldots,u_d)\in \R^d$.
Then
\begin{equation}\label{eq:bili}
\vv^TA_t\vu=\sum_{j+k=t+1}v_ju_k=0,\quad t=1,\ldots,d
\end{equation}
 implies that $v_1=u_1=0$.

Hence, according to the result above, the set
\[
\{\vv+\vu\in \R^d: (\vv,\vu) \text{ satisfies (\ref{eq:bili})}, \vv\neq 0, \vu\neq 0\}
\]
has zero Lebesgue outer measure in $\R^{d}$.  Then Lemma \ref{prop-PositiveMeasure} implies the $\A=(A_t)_{t=1}^d$ has almost everywhere phase retrieval property.
\eproof

\medskip

We next turn our attention to complex phase retrieval.

\medskip
\noindent
{\bf Proof of Theorem \ref{theo-ComplexGenericPRae}:}~~As before the random case is a corollary of the results on generic measurements, so we only need to prove the theorem for the cases of generic measurements.

First we consider the case where $V_j$ is the set of rank $r_j$ Hermitian matrices. Define the linear map $\varphi: \C^{d\times d} \lra \C^{d\times d}$ by
\begin{equation}  \label{HermitianIsom}
    \varphi(A) = \frac{1}{2} (A+A^T) +\frac{i}{2} (A-A^T).
\end{equation}
Then $\varphi$ is an isomorphism on $\C^{d\times d}$ with $\varphi^{-1}(B) =\frac{1}{2} (B+B^T)-\frac{i}{2}(B-B^T)$, and furthermore $\varphi$ restricted on $\R^{d\times d}$ is an isomorphism from $\R^{d\times d}$ to $\Herm(\C)$. For any $A,Q \in \C^{d\times d}$ define $L(A,Q):=\tr(\varphi(A)Q)$.

For any $s\geq 1$, let $V_s$ denote the set of matrices $A \in \C^{d\times d}$ such that $\rank(\varphi(A))\leq s$. The $V_s$ is clearly a projective variety. It was shown in \cite{WaXu16} that $V_s$ is admissible with respect to
$$
\{f^Q(A):=\tr(\varphi(A)Q):~0\neq Q\in {\mathcal M}_{d,2}(\C)\}
$$
for any $1 \leq s \leq d$. Let $V:=V_{r_1}\times \cdots\times V_{r_N}$. Thus, by Theorem \ref{theo-admissible} there exists a proper subvariety $Z$ of $V$ such that for any $\A=(A_j)_{j=1}^N \in V\setminus Z$ the projective variety
$$
      X_{\A}  := \Bigl\{Q\in {\mathcal M}_{d,2}(\C):~\tr(\varphi(A_j)Q) = 0 \mbox{~for all $1\leq j \leq N$}\Bigr\}
$$
has dimension $\dim(X_\A) = 4d-4-N\leq 2d-4$. Set $\varphi(\A):= (\varphi(A_j))_{j=1}^N$. In particular, if $\A$ is real and hence $\varphi(\A)\in \Herm^N(\C)$, by Theorem \ref{theo-PRae-VarietyC}, $\varphi(\A)$ has the almost everywhere PR property in $\C^d$.

Since $\varphi$ is an isomorphism on $\C^{d\times d}$, we have $\dim V_s = 2ds-s^2$. Moreover, $\dim_\R(V_s\cap\R^{d\times d})$ is exactly the (real) dimension of the set of Hermitian matrices with rank no more than $ s$, which is also $2ds-s^2$ (see also \cite[Lemma II.1]{KW15}). Thus,
$$
   \dim_\R(Z\cap(\R^{d\times d})^N)\leq \dim(Z) <\dim(V)= \dim_\R(V\cap (\R^{d\times d})^N).
$$
 For any $\A=(A_j)_{j=1}^N \in (V\setminus Z)\cap (\R^{d\times d})^N$ we have $\dim(X_\A) = 4d-4-N\leq 2d-4$. It follows from Theorem \ref{theo-PRae-VarietyC} that a generic  $\varphi(\A) = (\varphi(A_j))_{j=1}^N\in \Herm^N(\C)$ has the almost everywhere phase retrieval property. This proves the almost everywhere PR property in $\C^d$.

We now prove the case for orthogonal projections with prescribed ranks. Here the proof is virtually identical to the previous case, so we shall be rather brief. Similar to before, let $\varphi: \C^{d\times d} \lra \C^{d\times d}$ be defined by (\ref{HermitianIsom}), and for any $A,Q \in \C^{d\times d}$ define $L(A,Q):=\tr(\varphi(A)Q)$.

From now on the proof is almost verbatim from the proof of the previous case. For any $s\geq 1$,
we slightly abuse the notation and let $V_s$ denote the set of matrices $A$ in $\C^{d\times d}$ such that $\varphi^2(A)=\varphi(A)$ and $\rank(\varphi(A))\leq s$. It was shown in \cite{WaXu16} that $\dim(V_s) = 2s(d-s)+1$ and $V_s$ is admissible with respect to $\{f^Q(A):=\tr(\varphi(A)Q):~0\neq Q\in {\mathcal M}_{d,2}(\C)\}$ for any $1 \leq s \leq d$. By Theorem \ref{theo-admissible} there exists a subvariety $Z$ of $V=V_{r_1}\times \cdots\times V_{r_N}$ with $\dim Z <\dim V$ such that for any $\A=(A_j)_{j=1}^N \in V\setminus Z$, the subvariety of ${\mathcal M}_{d,2}(\C)$
$$
      X_{\A}  := \Bigl\{Q\in {\mathcal M}_{d,2}(\C):~\tr(\varphi(A_j)Q) = 0 \mbox{~for all $1\leq j \leq N$}\Bigr\}
$$
has dimension $\dim(X_\A) = 4d-4-N\leq 2d-4$.

As before, note that $\varphi(V_s\cap\R^{d\times d})$ is precisely the set of (complex) orthogonal projection matrices of rank $s$, which also has real dimension $2s(d-s)+1$.  Thus,
\[
\dim_\R(Z\cap(\R^{d\times d})^N)\leq \dim(Z)<\dim(V) \leq  \dim_\R(V\cap (\R^{d\times d})^N).
 \]
 For any $\A=(A_j)_{j=1}^N \in (V\setminus Z)\cap (\R^{d\times d})^N)$ we have $\dim(X_\A) = 4d-4-N\leq 2d-4$, and hence $\dim(\varphi(X_\A))\leq 2d-4$. It follows from Theorem \ref{theo-PRae-VarietyC} that $\varphi(\A) = (\varphi(A_j))_{j=1}^N$ has the almost everywhere phase retrieval property. We arrive at the conclusion.
\eproof

Theorem \ref{theo-ComplexGenericPRae} shows that $N\ge 2d$ generic Hermitian matrices or orthogonal projection matrices have almost everywhere phase retrieval property in $\C^d$. Note that Corollary \ref{coro-jacobi} states that $2d-1$ matrices are necessary to guarantee PR-ae property in $\C^d$. Then one may be interested in  whether there exist  $2d-1$ matrices which have PR-ae property in $\C^d$ ? Theorem \ref{th:thd} shows it is possible, which implies that the bound $2d-1$ is sharp.

\medskip
\noindent

{\bf Proof of Theorem \ref{th:thd}:}~~
To state conveniently, we use $\ve_j, j=1,\ldots,d,$ to denote the $d$-dimensional vector with the $j$-th entry being $1$ and other entries being $0$.
We construct $2d-1$ measurement matrices $A_j$ as follows:
\[
A_1=\ve_1\ve_1^T,\,\, A_j=\ve_1\ve_j^T+\ve_j\ve_1^T,\,\, A_{d-1+j}=i\ve_1\ve_j^T-i\ve_j\ve_1^T,\,\, j=2,\ldots,d.
\]
A simple observation is that $A_1,\ldots,A_{2d-1} \in \Herm(\C)$. For any $\vv=(v_1,\ldots,v_d) ,\vu=(u_1,\ldots,u_d)\in \C^d$ and $\vv\neq 0, \vu\neq 0$,  the equations
\begin{equation} \label{eq:C(2d-1)}
\Re(\vv^*A_j\vu)=0,\quad j=1,\ldots,2d-1,
\end{equation}
gives  that
\begin{equation}\label{eq:uiv}
\begin{aligned}
  \Re(\bar{v}_1u_1)&=0,   \\
     \Re(\bar{v}_1u_k+\bar{v}_ku_1)&=0,\quad  k=2,\ldots,d. \\
     \Re(i\bar{v}_1u_k-i\bar{v}_ku_1)&=0,
\end{aligned}
\end{equation}
We claim that if $u_1\neq 0$ and $v_1\neq 0$, then the solution to (\ref{eq:uiv}) satisfies
$\vu=ic\vv$ where $c\in \R$.  Based on Lemma \ref{prop-PositiveMeasure}, we just need to consider the case where either $u_1=0$ or $v_1=0$.
First, for the equation $\Re(\bar{v}_1u_1)=0$, if $u_1=0$, then the rest of equations gives $\bar{v}_1u_k=0$ for all $k=2,\ldots,d$, which implies $v_1=0$ since $\vu\neq 0$.
Similarly, if $v_1=0$, we can obtain $u_1=0$.
According to Lemma \ref{prop-PositiveMeasure}, if $v_1=u_1=0$ then ${\mathcal A}$ has the almost phase retrieval property. We arrive at the conclusion.

 We still need to prove $\vu=ic\vv$ if $u_1\neq 0$ and $v_1\neq 0$. Assume that $u_1\neq 0$ and $v_1\neq 0$. Then the equation $\Re(\bar{v}_1u_1)=0$ is equivalent to $u_1=icv_1$ for some real number $c\neq 0$, which
implies that $u_{1R}=-cv_{1I}$ and $u_{1I}=cv_{1R}$. Here, we denote $u_{1R}=\Re(u_1)$, $u_{1I}=\Im(u_1)$ and define $v_{1R}$, $v_{1I}$ similarly.

Next, for  $k=2,\ldots,d $, we consider the following two equations
\[
\Re(\bar{v}_1u_k+\bar{v}_ku_1)=0
\]
and
\[
\Re(i\bar{v}_1u_k-i\bar{v}_ku_1)=0.
\]
Putting $u_{1R}=-cv_{1I}$ and $u_{1I}=cv_{1R}$ into the above equations, we can obtain
\[ \left[\begin{array}{cc}
      v_{1R} & v_{1I} \\
      v_{1I} & -v_{1R}
    \end{array}\right]\left[\begin{array}{c}
                   u_{kR} \\
                   u_{kI}
                 \end{array}\right]=\left[\begin{array}{c}
                                 cv_{kR}v_{1I}-cv_{kI}v_{1R} \\
                                 -cv_{kR}v_{1R}-cv_{kI}v_{1I}
                               \end{array}
                 \right].
\]
Since $v_1\neq 0$, the above equations have a unique solution $u_{kR}=-cv_{kI},\; u_{kI}=cv_{kR}$. It implies that $u_k=icv_k$ for all $k=2,\ldots,d $. Thus, it gives that $\vu=ic\vv $ if $u_1\neq 0$ and $v_1\neq 0$.
\eproof

\section{Additional Results for  Almost Everywhere Standard Phase Retrieval}
\setcounter{equation}{0}

Finally we go back to the standard phase retrieval setting to tie up some loose ends for almost everywhere phase retrieval. Recall that in the standard setting each $A_j$ is a rank one matrix of the form $A_j = \vf_j\vf_j^*$, where $\vf_j\in\H^d$.   We say a group of vectors $\vf_1, \vf_2, \dots, \vf_{N}$ in $\H^d$ have the {\em almost everywhere phase retrieval property (PR-ae)} if and only if $\A=(A_j)_{j=1}^N$ has the property. Note that we often identify  $\vf_1, \vf_2, \dots, \vf_{N}$ with the $d\times N$ matrix $F = (\vf_1, \vf_2, \cdots, \vf_{N})$. Thus, we shall say $F$ has the almost everywhere phase retrieval property if $\vf_1, \vf_2, \dots, \vf_{N}$ have the property. In the real case, the PR-ae property has a simple characterization.

\begin{theo} \label{standPRae-real}
     Let $\vf_1, \vf_2, \dots, \vf_{N}$ be vectors in $\R^d$. For any $J \subset \{1, 2, \dots, N\}$ let $V_J$ denote the subspace $\Span \{\vf_j: j\in J\}$. Then $\vf_1, \vf_2, \dots, \vf_{N}$ have the PR-ae property if and only if for any $I, J \subset \{1, 2, \dots, N\}$ with $I \cup J = \{1,2,\dots, N\}$,
$$
       V_I^\perp + V_J^\perp \neq \R^d.
$$
In particular, $\vf_1, \vf_2, \dots, \vf_{d+1}\in \R^d$ have the PR-ae property if and only if they are full spark, i.e. any $d$ vectors among them are linearly independent.
\end{theo}
\Proof ($\Rightarrow$)~~Let $\vf_1, \vf_2, \dots, \vf_{N}$ have the PR-ae property. Assume that
$$
       V_I^\perp + V_J^\perp = \R^d.
$$
for some $I,J$ with $I \cup J = \{1,2, \dots, N\}$. We derive a contradiction. For any $\vu\in\R^d$ we have $\vu = \vx+\vy$ where $\vx\in V_I^\perp$ and $\vy \in V_J^\perp$. Set $\vv = \vx-\vy$. Then $\inner{\vf_i, \vu} = \pm \inner{\vf_i, \vv}$, depending on $i \in J$ or $i\in I$. Thus, $|\innerp{\vf_i, \vu}| = |\innerp{\vf_i, \vv}|$ for all $i$. This contradicts the PR-ae property of $\{\vf_i\}$.

($\Leftarrow$)~~Conversely, assume that $V_I^\perp + V_J^\perp \neq \R^d$ for all $I,J$ with $I \cup J = \{1,2, \dots, N\}$. We show that $\{\vf_i\}$ has the PR-ae property. Assume not, there exists a set $\Omega\subseteq \R^d$ with positive Lebesgue measure such that for each $\vu\in\Omega$ there exists a $\vv_\vu \neq \vu$ such that $\innerp{\vf_i, \vu} = \pm \innerp{\vf_i, \vv_\vu}$ for all $i$. Let
$$
    I_\vu:= \Bigl\{i:~\innerp{\vf_i, \vu} = -\innerp{\vf_i, \vv_\vu}\Bigr\}, \mhsp
    J_\vu:= \Bigl\{i:~\innerp{\vf_i, \vu} = \innerp{\vf_i, \vv_\vu}\Bigr\}.
$$
Clearly $I_\vu \cup J_\vu = \{1, \dots, N\}$. Since there are only finitely many distinct $I_\vu$ and $J_\vu$ there exist $I, J$ with $I \cup J = \{1, 2\dots, N\}$ and an $\tilde\Omega\subseteq \Omega$ of positive Lebesgue measure such that $I_\vu=I$ and $J_\vu = J$ for all $\vu\in\tilde\Omega$. Now for each $\vu\in \tilde\Omega$ we have $\vu+\vv_\vu\in V_I^\perp$ and $\vu-\vv_\vu\in V_J^\perp$. It follows that
$$
   \vu = \frac{1}{2}(\vu+\vv_\vu) + \frac{1}{2}(\vu-\vv_\vu) \in V_I^\perp + V_J^\perp.
$$
However,  $V_I^\perp + V_J^\perp$ is a subspace of $\R^d$ which contains a positive measure subset $\tilde\Omega$. Thus, it must be the whole space, namely $V_I^\perp + V_J^\perp = \R^d$. This is a contradiction.

Finally, we show if $\vf_1, \vf_2, \dots, \vf_{d+1}\in \R^d$ are full spark then for any $I, J \subset \{1, 2, \dots, d+1\}$ with $I \cup J = \{1,2,\dots, d+1\}$,
$$
       V_I^\perp + V_J^\perp \neq \R^d.
$$
First it is obvious that the conclusion holds provided $I=\{1, 2, \dots, d+1\}$ or $J=\{1, 2, \dots, d+1\}$.
We next consider the case where $I, J \subset \{1, 2, \dots, d\}$. Denote $n_1:=\dim(V_I)$ and $n_2:=\dim(V_J)$.
Since $I \cup J = \{1,2,\dots, d+1\}$ and $\vf_1, \vf_2, \dots, \vf_{d+1}$ are full spark, it implies that $n_1+n_2 \ge d+1$. Thus,
\[
\dim(V_I^\perp + V_J^\perp)=(d-n_1)+(d-n_2)\le d-1,
\]
which means that $ V_I^\perp + V_J^\perp \neq \R^d$. We arrive at the conclusion.
\eproof

We now turn our attention to the complex case. Theorem \ref{theo-ComplexGenericPRae}  shows that for $N$ generic $\vf_1, \vf_2, \dots, \vf_N$ where $N\geq 2d$,  $\A=(\vf_j\vf_j^*)_{j=1}^N$ has the PR-ae property. But is $N=2d$ minimal, namely, whether it is true that for $N \leq 2d-1$ we can never get the PR-ae property in the standard phase retrieval setting? This question was addressed in \cite{FSC05}. However, there was a general feeling on the mathematics perspective that the proof was not rigorous, and no one seemed to be able to verify its correctness, see Mixon \cite{Mix14}. Thus, this problem is still viewed, at least on the mathematics sideperspective, as an open problem. In this section we prove that a generic set $\{\vf_1, \vf_2, \dots, \vf_{2d-1}\}$ in $\C^d$ does not have the PR-ae property.

Note that, by Lemma \ref{prop-PositiveMeasure}, $\{\vf_1, \vf_2, \dots, \vf_{2d-1}\}$ in $\C^d$ has the PR-ae property if and only if the following holds: Let ${\mathcal U}$ be the set of $(\vu,\vv)\in \C^d\times\C^d$
such that $\vu \neq ic\vv$ for $c\in\R$ and
\begin{equation*}
       \Re (\overline{\inner{\vv, \vf_j}}\inner{\vu, \vf_j})=0 \mhsp
       \mbox{for all $1\le j \le 2d-1$.}
\end{equation*}
Then the set $E=\{\vu+\vv:~(\vu,\vv)\in {\mathcal U}, \vu\neq 0, \vv\neq 0\}$ has positive measure in $\C^d$.

To prove that a generic frame $\{\vf_1, \vf_2, \dots, \vf_{2d-1}\}$ in $\C^d$ does not have the almost everywhere phase retrieval property, we first consider the special case where $\vf_j = {\mathbf e}_j$ for $1 \leq j \leq d$. Thus, the associate frame matrix has the form
\begin{equation}  \label{frame-canonical}
      F=[I_d, G], \mhsp G \in\C^{d\times (d-1)}
\end{equation}
where $I_d$ is $d\times d$ identity matrix.
We shall prove that for a generic $G$ the frame is not PR-ae.

\begin{lem}  \label{lem-Remove1stRow}
   Assume that $F=[I_d,G]$ is a PR-ae in $\C^d$ where $G \in \C^{d\times (d-1)}$.
   Then $F_1=[I_{d-1}, G_1]$ is not a PR-ae in $\C^{d-1}$, where
   $G_1$ is obtained by removing the first row of $G$ (i.e., $F_1$ is $F$ with the first row being removed).
\end{lem}
\Proof  For any $\vu\in\C^d$, let $V_\vu$ denote the set of all $\vv\in\C^d$ such that
\begin{equation}  \label{2.5}
       \Re (\overline{\inner{\vu, \vf_j}}\inner{\vv, \vf_j})=0, \mhsp
       j=1,2, \dots, 2d-1
\end{equation}
where $\vf_j$ are the columns of $F$.
We shall only consider those $\vu$ with $u_1=0$, i.e. $\vu\in\{0\}\times\C^{d-1}$. Thus, for each such $\vu$ there exists a nonzero $\vv$ satisfying (\ref{2.5}). To see that, we only need to observe $\vv$ has $2d$ (real) unknowns while
there are only $2d-2$ linear equations in (\ref{2.5})  (we get $j=1$ for free because $u_1=0$).
Hence, $\dim_\R V_\vu \geq 2$. Let
$$
    Y = \left\{\vu:~V_\vu\subset\{0\}\times\C^{d-1}\right\}.
$$
We next divide the rest of the proof into two cases.

\smallskip
\noindent
{\bf Case I.}~~{\em $Y$ is not a null set in $\{0\}\times\C^{d-1}$.}

Note that if $\vv\in V_\vu$ is colinear with $\vu$ then we must have
$\vv=ic\vu$ for some $c\in\R$. Since $\dim_\R(V_\vu)\geq 2$ there must exists at least
a $\vv\in V_\vu$ that is not colinear with $\vu$. Now for any $\vu\in Y$,
since both $\vu,\vv$ have first entry 0, we have
\begin{equation}  \label{5.2}
       \Re (\overline{\inner{\vu', \vf_j'}}\inner{\vv', \vf_j'})=0 \mhsp
       \mbox{for all $j\geq 2$}
\end{equation}
where $\vu',\vv'$ are obtained by removing the first entry of $\vu,\vv$ and $\vf_j'$ are the
columns of $F_1$. However,  the set $\{\vu'\}$ satisfying (\ref{5.2}) has positive Lebesgue measure in $\C^{d-1}$, thus, $\{\vf_j'\}$ cannot have the PR-ae property. To see this more clearly, set $\vx'=\vu'+\vv'$ and $\vy'=\vu'-\vv'$. Then $|\inner{\vx', \vf_j'}| = |\inner{\vy', \vf_j'}| $. Note that at least one of the sets $\{\vx'\}$ and $\{\vy'\}$ has positive Lebesgue measure. We arrive at the conclusion in this case.

\smallskip
\noindent
{\bf Case II.}~~{\em $Y$ is a null set in $\{0\}\times\C^{d-1}$.}

We shall show that this case is impossible.
Without loss of generality we may simply assume that
$Y=\emptyset$. Thus, for any $\vu \in\{0\}\times\C^{d-1}$ there exists some
$\vv\in V_\vu$ with $v_1\neq 0$.
Clearly such a $\vv$ cannot be colinear with $\vu$ in $\C^d$.
Pick one such vector and denote it by $\vv_\vu$. Observe that $i\vu\in V_\vu$.
Hence, $V_\vu$ contains $t\vv_\vu +is\vu$ for all $t,s\in\R$.

We prove that $E:=\{\vu+\vv:~\vu\in \{0\}\times\C^{d-1}, \vv\in V_\vu\}$
has positive Lebesgue measure in $\C^d$, which implies that $F$ is not PR-ae in $\C^d$ (see Lemma \ref{prop-PositiveMeasure}). So, we have  a contradiction.
Note that for any nonzero $c\in\C$ we have $V_{c\vu} = cV_\vu$. Hence,
$V_{c\vu}$ contains $ct\vv_\vu +ics\vu$ for all $t,s\in\R$. Therefore, for each
$\vu\in \{0\}\times\C^{d-1}$ the set $E$ contains
$$
    c\vu + c(t\vv_\vu +ics\vu) = c(1+is)\vu + ct\vv_\vu
$$
for all $c\in\C$ and $t,s\in\R$. Set $c= (1+is)^{-1}$. It follows that
$\vu+ t(1+is)^{-1}\vv_\vu \in E$ for all $t,s\in\R$. But $z=t(1+is)^{-1}$ can take
on any non-imaginary complex number. Hence, $\vu+z\vv_\vu\in E$ for
any $z$ with $\Re(z) \neq 0$. In other words,
$$
   E \supseteq \left\{\vu+z\vv_\vu:~\vu\in \{0\}\times\C^{d-1}, z\in\C, \Re(z) \neq 0\right\}.
$$
Since each $\vv_\vu$ has nonzero first entry, it is clear that $E$ must have
positive Lebesgue measure in $\C^d$. This is a contradiction.
\eproof


%

\vspace{2mm}

\begin{theo} \label{theo-5.3}
A generic $F\in \C^{d\times (2d-1)}$ is not PR-ae.
\end{theo}
\Proof  It is well known (\cite{BCE06}) that the PR-ae property is preserved under nonsingular affine transformations. In other words, for any nonsingular $B\in \C^{d\times d}$, the frame $\{f_j\}_{j=1}^N$ has the PR-ae property if and only if $\{Bf_j\}_{j=1}^N$ does. Now for any nonsingular   $B$ and a generic $G\in \C^{d\times (d-1)}$ the frame $F = B \, [I_d, G]=[B, BG]$ is not PR-ae (otherwise,
Lemma \ref{lem-Remove1stRow} implies that a generic $F_1\in\C^{(d-1)\times (2(d-1)) }$ is not PR-ae which contradicts with Theorem \ref{theo-ComplexGenericPRae}
). Thus, a generic $F\in \C^{d\times (2d-1)}$ is not PR-ae.
\eproof

\end{document}